\documentclass[12pt, reqno]{amsart}
\usepackage{amsmath, amsthm, amscd, amsfonts, amssymb, graphicx, color}
\usepackage[bookmarksnumbered, colorlinks, plainpages]{hyperref}

\newtheorem{theorem}{Theorem}[section]
\newtheorem{lemma}[theorem]{Lemma}

\theoremstyle{definition}
\newtheorem{definition}[theorem]{Definition}
\newtheorem{example}[theorem]{Example}

\theoremstyle{remark}

\numberwithin{equation}{section}

\input{mathrsfs.sty}

\begin{document}
\setcounter{page}{1}

\title[Amenable group and approximate identities]{Amenable groups and bounded $\Delta$-weak approximate identities}

\author[M.fozouni]{Mohammad Fozouni}

\address{Faculty of Mathematics and computer science, Department of mathematics, Kharazmi University,
599 Taleghani, Tehran, Iran.}
\email{\textcolor[rgb]{0.00,0.00,0.84}{fozouni@khu.ac.ir or fozouni@hotmail.com}}

\subjclass[2010]{22D05, 43A15, 43A07.}

\keywords{Amenable group, Character, Figa-Talamanca Herz algebra.}

\begin{abstract}
Let $A$ be a Banach algebra with a non-empty character space. We say that a bounded net $\{e_{\alpha}\}$ in $A$ is a bounded $\Delta$-weak approximate identity for $A$ if, for each $a\in A$ and compact subset $K$ of $\Delta(A)$, $||\widehat{e_{\alpha}a}-\widehat{a}||_{K}=\sup_{\phi\in K}|\phi(e_{\alpha}a)-\phi(a)|\rightarrow 0$. For each $1<p<\infty$, we prove that the Figa-Talamanca Herz algebra, $A_{p}(G)$ has a bounded $\Delta$-weak approximate identity if and only if $G$ is an amenable group. Also we give a sufficient condition for amenability of group $G$.
\end{abstract} \maketitle

\section{Introduction and preliminaries}

\noindent Let $A$ be a Banach algebra and let $\Delta(A)$ be the character space of $A$, i.e., the space consisting of all nonzero homomorphisms from $A$ into $\mathbb{C}$. Throughout this paper we will always let $A$ be a Banach algebra which $\Delta(A)$ is nonempty.

For the convenience of the reader here we give a brief outline of important definitions and results.

Let $G$ be a locally compact group. If $1<p<\infty$, we let $A_{p}(G)$ denote the subspace of $C_{0}(G)$ consisting of functions of the form, $u=\sum_{i=1}^{\infty}f_{i}\ast \widetilde{g_{i}}$, where $f_{i}\in L^{p}(G)$, $g_{i}\in L^{q}(G)$, $1/p+1/q=1$, $\sum_{i=1}^{\infty}||f_{i}||_{p}||g_{i}||_{q}<\infty$ and $\widetilde{f}(x)=\overline{f(x^{-1})}$ for all $x\in G$. $A_{p}(G)$ is called the Figa-Talamanca Herz algebra and with the pointwise operation and the following norm is a Banach algebra.
\begin{equation*}
||u||_{A_{p}(G)}=\inf\{\sum_{i=1}^{\infty}||f_{i}||_{p}||g_{i}||_{q}: u=\sum_{i=1}^{\infty}f_{i}\ast \widetilde{g_{i}}\}\cdot
\end{equation*}
It is obvious that for each $u\in A_{p}(G)$, $||u||\leq ||u||_{A_{p}(G)}$, where $||u||$ is the norm of $u$ in $C_{0}(G)$. Also we know that $\Delta(A_{p}(G))=G$, i.e., each character of $A_{p}(G)$ is an evaluation function at some $x\in G$ \cite[Theorem 3]{Herz}.

The group $G$ is said to be amenable if there exists an $m\in L^{\infty}(G)^{*}$ such that $m\geq 0$, $m(1)=1$ and $m(L_{x}f)=m(f)$ for each $x\in G$ and $f\in L^{\infty}(G)$, where $L_{x}f(y)=f(x^{-1}y)$ \cite[Definition 4.2]{Pier}.

There are many characterizations of amenability of a group $G$ that can be found in the literature. One of these characterizations is the following \cite[Theorem 9.6]{Pier}. Here $C_{c}^{+}(G)$ denote the space of all continuous functions from $G$ into $\mathbb{C}$ such that   are positive and have compact support.
\begin{theorem}\label{maintheorem} Let $G$ be a locally compact group. The group $G$ is amenable if and only if for one $q\in (1,\infty)$ and each $f\in C_{c}^{+}(G)$
\begin{align*}
||f||_{1}=||L_{f}||_{CV_{q}(G)}&=\sup\{||L_{f}(g)||:g\in L^{q}(G), ||g||_{q}\leq 1\}\\
&=\sup\{||f\ast g||:g\in L^{q}(G), ||g||_{q}\leq 1\}\cdot
\end{align*}
\end{theorem}
It is easy to verify that, $C_{c}(G)\ast C_{c}(G)\subseteq C_{c}(G)$ and for each $n\in \mathbb{N}$ and $\psi\in C_{c}^{+}(G)$, $||\psi||_{1}^{n}=||(\psi^{*})^{n}||_{1}$, which $\psi^{*}(x)=\Lambda(x^{-1})\overline{\psi(x^{-1})}$, such that $\Lambda$ show the modular function of group $G$.

If $\phi\in C_{c}(G)\subseteq M(G)$, the function $F_{\phi}:A_{p}(G)\rightarrow \mathbb{C}$, defined by,
 \begin{center}
 $F_{\phi}(u)=<u,\phi>\hspace{0.5cm}(u\in A_{p}(G)),$
 \end{center}
 is an element of $A_{p}(G)^{*}$. Here $<.,.>$ denote the pairing between $M(G)$ and $C_{0}(G)$. In view of \cite[Proposition 10.3]{Pier} we have $||F_{\phi}||=||L_{\phi}||_{CV_{q}(G)}$.

It is easy to see that for each $\psi\in C_{c}^{+}(G)$, $n\in \mathbb{N}$ and  $p\in(1,\infty)$,
\begin{equation*}
||L_{(\psi^{*})^{n}}||_{CV_{p}(G)}\leq ||L_{\psi}||_{CV_{p}(G)}^{n}\cdot
\end{equation*}
Let $A$ be a commutative Banach algebra. A net $\{u_{\alpha}\}$ in $A$ is called a bounded weak approximate identity, if there exists a nonnegative constant $C<\infty$ such that $||u_{\alpha}||<C$ for each $\alpha$ and
\begin{equation*}
\lim_{\alpha}|\phi(au_{\alpha})-\phi(a)|=0,
\end{equation*}
 for all $a\in A$ and $\phi\in\Delta(A)$.

 Weak approximate identities were introduced by C.Jones and C.Lahr in \cite{Jones}. The authors give an example of a Banach algebra that has a bounded weak approximate identity, but does not have an approximate identity, bounded or unbounded.

In a similar fashion we define the concept of $\Delta$-weak approximate identity for arbitrary Banach algebras (not necessarily commutative), but no attempt has been made here to show the difference between these two kinds of approximate identities. In the following section, first we give the definition of a bounded $\Delta$-weak approximate identity for a Banach algebra. Then we give an examples of  a Banach algebra such that has  bounded $\Delta$-weak approximate identity, but does not have any bounded approximate identity. We show that for Figa-Talamanca Hera algebras, the existence of this kind of approximate identity yields it has a bounded approximate identity also. We prove this, with use of the structure of amenable groups. Finally we prove that if for some $\varphi\in \Delta(A_{p}(G))\cup\{0\}$, $\ker(\varphi)$ has a bounded $\Delta$-weak approximate identity, then $G$ is an amenable group.

\section{Main Results}

We commence this section with the definition of the $\Delta$-weak approximate identity for a Banach algebra $A$ as follows. Note that for each $a\in A$, $\widehat{a}$, denote the Gelfand's transform of $a$.

\begin{definition}\label{maindef}
Let $A$ be a Banach algebra. A  \emph{$\Delta$-weak approximate identity} for $A$ is a net $\{e_{\alpha}\}$ such that for each $a\in A$ and compact subset $K$ of $\Delta(A)$
\begin{equation*}\label{formula1}
\lim_{\alpha}||\widehat{ae_{\alpha}}-\widehat{a}||_{K}=\lim_{\alpha}\sup_{\phi\in K}|\phi(ae_{\alpha})-\phi(a)|=0\cdot
\end{equation*}
If the net $\{e_{\alpha}\}$ is bounded, we say that it is a bounded  $\Delta$-weak approximate identity (b.$\Delta$-w.a.i) for Banach algebra $A$.

For a subset $B$ of $A$ we say $B$ has a $\Delta$-w.a.i if there exists a net $(b_{\alpha})$ in $B$ such that for each $b\in B$ and compact subset $K$ of $\Delta(A)$
\begin{equation*}
\lim_{\alpha}||\widehat{bb_{\alpha}}-\widehat{b}||_{K}=0\cdot
\end{equation*}

\end{definition}
Here the definition of one sided (left/right) $\Delta$-weak approximate identity are redundant, because each $\phi\in\Delta(A)$ is multiplicative.

Clearly, if a Banach algebra $A$  has a left(right or two sided ) identity, then it has a bounded $\Delta$-weak approximate identity. Because, if $e$ is a right identity of $A$, then consider the bounded net $\{e_{\alpha}\}$ with $e=e_{\alpha}$ for each $\alpha$. Hence for all $a\in A$ and $K\subseteq \Delta(A)$ we have
\begin{equation*}
||\widehat{ae_{\alpha}}-\widehat{a}||_{K}=\sup_{\phi\in K}|\phi(ae_{\alpha})-\phi(a)|=0\\
\end{equation*}
The following example give a Banach algebra such that has a left identity, but does not have any approximate identity.
\begin{example} Let $X$ be a Banach space and $B=X^{*}$ be its topological dual. Take $a\in X$ with $||a||\leq 1$ and for each $f, g\in B$, define $f.g:=f(a)g$. It is obvious that $B$ with this multiplication is a Banach algebra, which we denote it by $B_{a}(X)$. The Banach algebra $B_{a}(X)$ has a bounded approximate identity if and only if $\textmd{dim}(X)=1$ \cite{Laali}. To see this, let $(b_{\alpha})$ be an  approximate identity for $B_{a}(X)$. So,
\begin{equation} \label{eq3}
||b_{\alpha}b-b||, ||bb_{\alpha}-b||\rightarrow 0\quad(b\in B_{a}(X)).
\end{equation}
Let $M=<a,a^{'}>$ be the subspace of $X$ generated by $a, a^{'}$. Then $b^{'}:M\longrightarrow \mathbb{C}$ defined by,
\begin{equation*}
b^{'}(\gamma_{1}a+\gamma_{2}a^{'})=\gamma_{2}||a^{'}||\quad(\gamma_{1}, \gamma_{2}\in \mathbb{C}),
\end{equation*}
is a bounded linear functional on $M$.  Therefore, the Hahn-Banach theorem implies that there exists a functional $b\in B$, such that, $b_{|M}=b^{'}$ and $||b||=||b^{'}||$. Hence $0\neq b\in B$ and $b(a)=0$. Therefore, Relation (\ref{eq3}) does not hold for any $b\in B$.

So, if we take $X$ a Banach space such that $\textmd{dim}(X)>1$, then $B_{a}(X)$ is a Banach algebra such that does not have any bounded approximate identity. Also each $b_{0}\in B_{a}(X)$ which $b_{0}(a)=1$ is a left identity for $B_{a}(X)$.  Also it is clear that $\Delta(B_{a})=\{\widehat{a}\}$.
\end{example}

\subsection{Main Theorems}

A classical theorem due to Leptin and Herz, characterize the amenability of a group $G$ through the existence of a bounded approximate identity for Figa-Talamanca Herz algebras as follows.
\begin{theorem}\label{T: Leptin-Herz}
Let $G$ be a locally compact group and for each $1<p<\infty$, $A_{p}(G)$ be the Figa-Talamanca Herz algebra. Then $A_{p}(G)$ has a bounded approximate identity if and only if $G$ is an amenable group.
\end{theorem}
The proof of the above theorem in case $p=2$ is due to Leptin \cite{Leptin} and in general is due to Herz \cite{Herz}.

Now we give the following theorem  that is one of our main results.
\begin{theorem}\label{T: main} Let $G$ be a locally compact group and $1<p<\infty$. Then the following are equivalent.
\begin{enumerate}
  \item $G$ is an amenable group,
  \item $A_{p}(G)$ has a b.$\Delta$-w.a.i,
  \item $A_{p}(G)$ has a b.a.i.
\end{enumerate}
\end{theorem}
\begin{proof} By use of Theorem \ref{T: Leptin-Herz} and this fact that each b.a.i for $A_{p}(G)$ is a b.$\Delta$-w.a.i, we have the implications $(1):\rightarrow (3):\rightarrow (2)$.

Now, let $(2)$ holds and let $\{u_{i}\}$ be a b.$\Delta$-w.a.i for $A_{p}(G)$ bounded by $C$. Suppose that $q$ is the conjugate exponent of $p$, i.e., $1/p+1/q=1$. By Theorem \ref{maintheorem}, it is enough to show that for $q$ and each $\psi\in C_{c}^{+}(G)$, $||L_{\psi}||_{CV^{q}}=||\psi||_{1}$.

If $K$ is any compact subset of $G$, we choose $V$ an arbitrary compact neighborhood of $G$ containing $e$ (identity of group) and put, $f=|V|^{-1}1_{V}\ast 1_{V^{-1}K}$.

 A routine verification shows that if $x\in K$, $f(x)=1$ and otherwise $f(x)=0$.
Since $\{u_{i}\}$ is a b.$\Delta$-w.a.i and $f\in A_{p}(G)$, for $K\subseteq G=\Delta(A_{p}(G))$, we have
\begin{equation*}
||\widehat{u_{i}f}-\hat{f}||_{K}=\sup_{t\in K}|u_{i}(t)f(t)-f(t)|=\sup_{t\in K}|u_{i}(t)-1|\rightarrow 0\cdot
\end{equation*}
Hence, for $\epsilon>0$, there exists $i_{0}$, such that, $\sup_{t\in K}|Re(u_{i_{0}}(t))-1|<\epsilon$. Therefore, $\inf\{Re(u_{i_{0}}(t));t\in K\}\geq 1-\epsilon$.

Let $\phi\in C_{c}^{+}(G)$ and $K=\textrm{supp}(\phi)$.
 By the discussion after Theorem \ref{maintheorem}, we have
\begin{equation*}
|<u_{i_{0}},\phi>|=|F_{\phi}(u_{i_{0}})|\leq ||L_{\phi}||_{CV_{q}(G)}||u_{i_{0}}||\leq C||L_{\phi}||_{CV_{q}(G)}\cdot
\end{equation*}
But
\begin{equation*}
Re<u_{i_{0}},\phi>=\int_{K} Re(u_{i_{0}}(x))\phi(x)dx\geq (1-\epsilon)||\phi||_{1}\cdot
\end{equation*}
Hence, if $\epsilon$ tend to 0, we have $||\phi||_{1}\leq C||L_{\phi}||_{CV_{q}(G)}$.

Let $\psi\in C_{c}^{+}(G)$ be arbitrary and $n\in \mathbb{N}$. Thus we have
\begin{equation*}
||\psi||_{1}^{n}=||(\psi^{*})^{n}||_{1}\leq C||L_{(\psi^{*})^{n}}||_{CV_{q}(G)}\leq C||L_{\psi}||_{CV_{q}(G)}^{n}\cdot
\end{equation*}
Therefore, $||\psi||_{1}\leq ||L_{\psi}||_{CV_{q}(G)}$.

The verification of $||L_{\psi}||_{CV_{q}(G)}\leq ||\psi||_{1}$ is easy. Hence, $||L_{\psi}||_{CV_{q}(G)}=||\psi||_{1}$ and this complete the proof.
\end{proof}
In the sequel we find a sufficient condition that in $A_{p}(G)$, $G$ is an amenable group.

First we mention the following lemma that is a variant of \cite[Proposition 33.2]{Doran}.
\begin{lemma}\label{L: main} Let $A$ be a Banach algebra and $I$ be a closed two-sided  ideal of $A$ which has a b.$\Delta$-w.a.i and the quotient Banach algebra $A/I$ has a  bounded left approximate identity (b.l.a.i). Then $A$ has a b.$\Delta$-w.a.i.
\end{lemma}
\begin{proof}
 Let $\{e_{\alpha}\}$ be a b.$\Delta$-w.a.i for $I$ and $\{f_{\delta}+I\}$ be a b.l.a.i for $A/I$. Suppose that $F=\{a_{1},...a_{m}\}$ is a finite subset of $A$ and $n$ is a positive integer. Let $M$ be an upper bound for $\{||e_{\alpha}||\}$. For $\lambda=(F,n)$, there exists $f_{\delta_{\lambda}}$ such that
\begin{equation*}
||f_{\delta_{\lambda}}a_{i}-a_{i}+I||<\frac{1}{2(1+M)n}\hspace{0.5cm}(i=1,2,3,...,m)\cdot
\end{equation*}
By properties of the quotient norm of $A/I$, there exists $y_{i}\in I$ such that
\begin{equation*}
||f_{\delta_{\lambda}}a_{i}-a_{i}+y_{i}||<\frac{1}{2(1+M)n}\hspace{0.5cm}(i=1,2,3,...,m)\cdot
\end{equation*}
Let $K\subseteq \Delta(A)$ be a compact set. Since $\{e_{\alpha}\}$ is a b.$\Delta$-w.a.i, for each $y_{i}$ with $i\in\{1,2,3,...,m\}$, which satisfy the above relation, there exists $e_{\alpha_{\lambda}}\in\{e_{\alpha}\}$ such that
\begin{equation*}
\sup_{\psi\in K}||\psi(e_{\alpha_{\lambda}}y_{i})-\psi(y_{i})||<\frac{1}{2n}\hspace{0.5cm}(i=1,2,3,...,m)\cdot
\end{equation*}
Now for each $i\in\{1,2,3,...,m\}$ we have
\begin{align*}
\sup_{\psi\in K}||\psi((e_{\alpha_{\lambda}}+f_{\delta_{\lambda}}-e_{\alpha_{\lambda}}f_{\delta_{\lambda}})a_{i})-\psi(a_{i})||&\leq \sup_{\psi\in K}||\psi(f_{\delta_{\lambda}}a_{i}-a_{i}+y_{i})||\\
&+\sup_{\psi\in K}||\psi(e_{\alpha_{\lambda}}y_{i})-\psi(y_{i})||\\
&+\sup_{\psi\in K}||\psi(e_{\alpha_{\lambda}}a_{i}-e_{\alpha_{\lambda}}f_{\delta_{\lambda}}a_{i}-e_{\alpha_{\lambda}}y_{i})||\\
&\leq ||f_{\delta_{\lambda}}a_{i}-a_{i}+y_{i}||\\
&+\frac{1}{2n}+M||a_{i}-f_{\delta_{\lambda}}a_{i}-y_{i}||\\
&<\frac{1}{n}\cdot
\end{align*}
Therefore $\{e_{\alpha_{\lambda}}+f_{\delta_{\lambda}}-e_{\alpha_{\lambda}}f_{\delta_{\lambda}}\}_{\lambda\in \Lambda}$ is a $\Delta$-w.a.i for $A$, that $\Lambda=\{(F,n):F\subseteq A \textmd{ is  finite}, n\in \mathbb{N}\}$ is a directed set with $(F_{1},n_{1})\leq (F_{2},n_{2})$ if, $F_{1}\subseteq F_{2}$ and $n_{1}\leq n_{2}$.

Now we show that there exists a b.$\Delta$-w.a.i for $A$. Since $\{f_{\delta}+I\}$ is bounded, there exists a positive integer $N$ such that, $||f_{\delta}+I||<N$ for each $\delta$. So in view of the properties of the quotient norm, there exists $y_{\delta}\in I$ such that, $||f_{\delta}+I||<||f_{\delta}+y_{\delta}||<N$. Put $f^{'}_{\delta}=f_{\delta}+y_{\delta}$. Hence $\{f^{'}_{\delta}+I\}$ is a bounded approximate identity for $A/I$ that $\{f^{'}_{\delta}\}$ is  bounded.
Now we have
\begin{center}
$||e_{\alpha_{\lambda}}+f^{'}_{\delta_{\lambda}}-e_{\alpha_{\lambda}}f^{'}_{\delta_{\lambda}}||\leq ||e_{\alpha_{\lambda}}||+||f^{'}_{\delta_{\lambda}}||+||e_{\alpha_{\lambda}}||||f^{'}_{\delta_{\lambda}}||<M+N+NM\cdot$
\end{center}
Therefore $A$ has a b.$\Delta$-w.a.i.
\end{proof}
We end this paper by the following theorem which provide a sufficient condition for amenability of group $G$ in  $A_{p}(G)$.
\begin{theorem} Let $G$ be a locally compact group, $1<p<\infty$ and $\varphi\in \Delta(A_{p}(G))\cup\{0\}$. If $\ker(\varphi)$ has a $b.\Delta.w.a.i$, then $G$ is an amenable group.
\end{theorem}
\begin{proof} Since $\ker(\varphi)$ is an ideal of $A_{p}(G)$ with codimension one, the quotient Banach algebra $A_{p}(G)/\ker(\varphi)$ has a b.a.i. So the result follows by use of Theorem \ref{T: main} and Lemma \ref{L: main}.
\end{proof}

\bibliographystyle{amsplain}

\end{document}